\documentclass[final,leqno,onefignum,onetabnum]{siamltex1213}

\usepackage[utf8]{inputenc}

\usepackage{inconsolata} 
\usepackage{listings}
\usepackage{booktabs}
\usepackage{algpseudocode}
\usepackage[table]{xcolor}

\definecolor{red2}{HTML}{FDBB84}
\definecolor{red1}{HTML}{E34A33}

\usepackage{fancyvrb}
\definecolor{objcol}{HTML}{CA0020}
\definecolor{flowconcol}{HTML}{F4A582}
\definecolor{totflowcol}{HTML}{92C5DE}
\definecolor{vardefcol}{HTML}{0571B0}

\lstdefinelanguage{julia}
{
  keywordsprefix=\@,
  morekeywords={
    exit,whos,edit,load,is,isa,isequal,typeof,tuple,ntuple,uid,hash,finalizer,convert,promote,
    subtype,typemin,typemax,realmin,realmax,sizeof,promote_type,method_exists,applicable,
    invoke,dlopen,dlsym,system,error,throw,assert,new,Inf,Nan,pi,im,begin,while,for,in,return,
    break,continue,macro,quote,let,if,elseif,else,try,catch,end,bitstype,ccall,do,using,module,
    import,export,importall,baremodule,immutable,local,global,const,Bool,Int,Int8,Int16,Int32,
    Int64,Uint,Uint8,Uint16,Uint32,Uint64,Float32,Float64,Complex64,Complex128,Any,Nothing,None,
    function,type,typealias,abstract
  },
  sensitive=true,
  morecomment=[l]{\#},
  morestring=[b]',
  morestring=[b]" 
}

\usepackage{amsmath}
\usepackage{amssymb}
\usepackage{nicefrac}
\usepackage{bm}
\usepackage{tikz,tikz-qtree}
\usepackage{multirow}
\usetikzlibrary{shapes}

\title{J\lowercase{u}MP: a modeling language for mathematical optimization} 

\author{Iain Dunning, Joey Huchette, Miles Lubin \thanks{MIT Operations Research Center \texttt{\{idunning,huchette,mlubin\}@mit.edu}.}}

\begin{document}
\maketitle
\slugger{sirev}{xxxx}{xx}{x}{x--x}

\begin{abstract}
JuMP is an open-source modeling language that allows users to express a wide range of optimization problems (linear, mixed-integer, quadratic, conic-quadratic, semidefinite, and nonlinear) in a high-level, algebraic syntax. JuMP takes advantage of advanced features of the Julia programming language to offer unique functionality while achieving performance on par with commercial modeling tools for standard tasks. In this work we will provide benchmarks, present the novel aspects of the implementation, and discuss how JuMP can be extended to new problem classes and composed with state-of-the-art tools for visualization and interactivity.
\end{abstract}



\pagestyle{myheadings}
\thispagestyle{plain}
\markboth{Dunning et al.}{JuMP}

\section{Introduction}

William Orchard-Hays, who developed some of the first software for linear programming (LP) in collaboration with George Dantzig, observed that the field of mathematical optimization developed hand-in-hand with the field of computing~\cite{OrH-History}.  
Beginning with the introduction of IBM's first commercial scientific computer in 1952, advancements in technology were immediately put to use for solving military and industrial planning problems.
LP software was viewed as generally reliable by the 1970s, 
when mainframe computers had become mainstream. 
However, developers of these systems recognized that the difficulty of translating the complex mathematical formulation of a problem into the requisite input formats based on punch cards was a major barrier to adoption~\cite{Fourerismp}.

In the late 1970s, the first algebraic modeling languages (AMLs) were developed with the aim of allowing users to express LP and other optimization problems in a natural, algebraic form similar to the original mathematical expressions, much in the same way that MATLAB was created contemporaneously to provide a high-level interface to linear algebra. 
Similar to how MATLAB translated user input into calls to LINPACK~\cite{linpack}, 
AMLs do not \textit{solve} optimization problems; they provide the problems to optimization routines called \textit{solvers}. 
GAMS~\cite{GAMS} and AMPL~\cite{AMPLBook}, two well-known commercial AMLs whose development started in 1978 and 1985 respectively, are widely recognized among similar systems like AIMMS, LINDO/LINGO, and MPL as having made a significant impact on the adoption of mathematical optimization in a number of fields.

In this paper, we present JuMP, an AML which is embedded in the Julia programming language~\cite{Julia}. 
In addition to providing a performant open-source alternative to commercial systems,
JuMP has come to deliver significant advances in modeling and extensibility by taking advantage of a number of features of Julia which are unique within the realm of programming languages for scientific computing. We highlight the novel technical aspects of JuMP's implementation in sufficient generality to apply broadly beyond the context of AMLs, in particular for the implementation of scientific domain-specific languages~\cite{dsl,STAN,UFL} and of automatic differentiation (AD) techniques for efficient computations of derivatives~\cite{Griewank2008EDP,NaumannAD}.

To date, AMPL, GAMS, and similar commercial packages represent the state of the art in AMLs and are widely used in both academia and industry. These AMLs
are quite efficient at what they were designed for;
however, a number of drawbacks motivated us to develop a new AML.
Unsatisfied with relatively standalone commercial systems, we wanted a lightweight AML which fits naturally within a modern scientific workflow. Such workflows could  include solving optimization problems within a larger simulation or interactive visualization, for example, or constructing a complex optimization model programmatically from modular components~\cite{Kallrath,gpkit}.
As algorithm developers, we wanted to be able to interact with solvers while they are running, for both control of the solution process and to reduce the overhead of regenerating a model when solving a sequence of related instances~\cite{GAMSGUSS}.
Finally, as modelers, we wanted to create user-friendly AML extensions for new problem classes that couple with specialized solution approaches; in contrast, commercial AMLs were not designed to be extended in this way~\cite{SET}. In short, with similar motivations as the developers of the Julia language itself~\cite{juliablog}, we created JuMP because we wanted more than what existing tools provided.

JuMP joins a rich family of open-source AMLs which been developed by academics since the 2000s. YALMIP~\cite{YALMIP} and CVX~\cite{cvx-matlab}, both based on MATLAB, were created to provide functionality such as handling of semidefinite and disciplined convex~\cite{DCP} optimization, which was not present in commercial AMLs. CVX in particular has been cited as making convex optimization as accessible from MATLAB as is linear algebra and was credited for its extensive use in both research and teaching~\cite{cvxprize}. 
Pyomo~\cite{Pyomo} is an AML which was originally designed to recreate the functionality of AMPL in Python and was later extended to new problem classes such as stochastic programming~\cite{PySP}.
Embedded within general-purpose programming languages\footnote{An idea which traces back to the commercial ILOG C++ interface of the 1990s}, these open-source AMLs broadly address our concerns of fitting within a modern workflow and are powerful and useful in many contexts.
However, their slow performance, due to being embedded in high-level languages like MATLAB and Python, motivated our preliminary work investigating Julia as an alternative high-level host language with the promise of fewer performance compromises~\cite{JuMPIJOC}.

Following JuMP's first release in 2013, which supported linear and mixed-integer optimization, we have enabled modeling for quadratic, conic-quadratic, semidefinite, and
general derivative-based nonlinear optimization problems, standard problem classes supported by the commercial and open-source AMLs. At the same time, we have extended the functionality of JuMP beyond what is typically available in an AML, either commercial or open-source. These features, which will be described in the text, include callbacks for in-memory bidirectional communication with branch-and-bound solvers, automatic differentiation of user-defined nonlinear functions, and easy-to-develop add-ons for specialized problem classes such as robust optimization. 
JuMP's unique mix of functionality has driven growing adoption by researchers~\cite{Frans,Giordano,nikita,vulnerability,Fusion}, and JuMP and has been used for teaching in courses in at least 10 universities (e.g.,~\cite{IAPClass}). In this paper, we will highlight both the important technical and usability aspects of JuMP, including how JuMP itself uses the advanced features of Julia.

The remainder of the paper is structured as follows. In Section~\ref{sec:role} we introduce in more detail the tasks required of an AML and present an example of AML syntax. In Sections \ref{sec:macros} and \ref{sec:linear} we discuss JuMP's use of syntactic macros and code generation, two advanced technical features of Julia which are key to JuMP's performance. In Section~\ref{section:nonlinear} we discuss JuMP's implementation of derivative computations. In Section~\ref{sec:extensions} we discuss a number of powerful extensions which have been built on top of JuMP, and in Section~\ref{sec:ijulia} we conclude with a demonstration of how JuMP can be composed with the growing ecosystem of Julia packages to produce a compelling interactive and visual user interface with applications in both academia and industry.

\section{The role of a modeling language}\label{sec:role}

\begin{figure}[t]
\centering
\fbox{\includegraphics[width=0.9\textwidth]{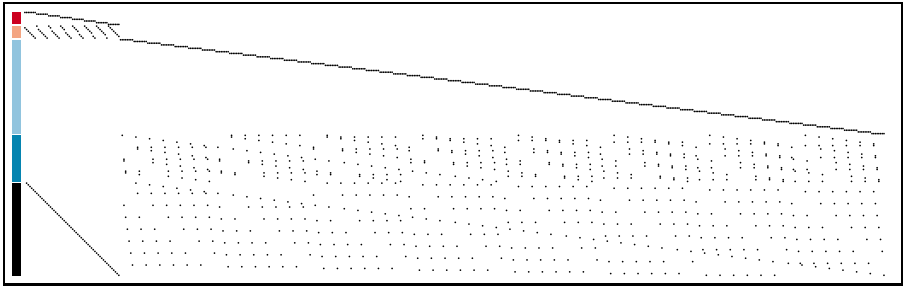}}
\caption{
Sparsity pattern from the constraint coefficient matrix for a  multicommodity flow problem arising from optimal routing in communication networks~\cite{Bienstock1995}. The dots correspond to nonzero elements of the matrix. We identify five groups of constraints, indicated with colored strips on the left. Modeling languages remove the need to write code to generate such complex matrices by hand; users instead work with a much more natural algebraic representation of the optimization problem.
}
\label{fig:sparsepattern}
\end{figure}

Prior to the introduction of AMLs (and continuing to a lesser degree today), users would write low-level code which directly generated the input data structures for an optimization problem.
Recall that standard-form linear programming problems can be stated as
\begin{equation}\label{LP-standard}
\begin{split}
\min_{x \in \mathbb{R}^n}\quad & c^Tx\\
\text{s.t.}\quad& Ax = b, \\
&x \geq 0,
\end{split}
\end{equation}
that is, minimization of a linear objective subject to linear equality and inequality constraints (all elements of $x$ must be nonnegative).
In the case of LP, the input data structures are the vectors $c$ and $b$ and the matrix $A$ in sparse format, and the routines to generate these data structures are called \textit{matrix generators}~\cite{Fourerismp}. Typical mathematical optimization models have complex indexing schemes; for example, an airline revenue management model may have decision variables $x_{s,d,c}$ which represent the number of tickets to sell from the source $s$ to destination $d$ in fare class $c$, where not all possible combinations of source, destination and fare class are valid. A matrix generator would need to efficiently map these variables into a single list of linear indices and then construct the corresponding sparse matrix $A$ as input to the solver, which is tedious, error-prone, and fragile with respect to changes in the mathematical model. An example sparsity pattern in Figure~\ref{fig:sparsepattern} demonstrates that these can be quite complex even for small problems. This discussion extends naturally to quadratic expressions $c^Tx + \frac{1}{2}x^TQx$; the matrix $Q$ is simply another component of the input data structure. The role of an AML in these cases is to accept closed-form algebraic expressions as user input and transparently generate the corresponding input matrices and vectors, removing the need to write matrix generators by hand. AMLs additionally handle any low-level details of communicating with the solver, either via a callable library or by exchanging specially formatted files.

AMLs have a similar role in the context of
nonlinear optimization problems, which often arise in scientific and engineering applications. The standard form for derivative-based nonlinear optimization problems is 
\begin{equation}\label{eq:nlp-standard}
\begin{split}
\min_{x \in \mathbb{R}^n}\quad & f(x)\\
\text{s.t. }\quad & g_i(x) \leq 0 \quad i = 1,\ldots,m_g,\\
&h_i(x) = 0  \quad i =1,\ldots, m_h,
\end{split}
\end{equation}
where $f, g_i, h_i : \mathbb{R}^n \to \mathbb{R}$ are linear or nonlinear functions. Depending on certain properties of $f$, $g$, and $h$ such as convexity, these problems may be easy or hard to solve to a global solution; regardless, the solution methods often rely on the availability of first-order derivatives, that is, the gradient vectors $\nabla f(x)$, $\nabla g_i(x)$, and $\nabla h_i(x)$, and may be further accelerated by the availability of second-order derivatives, that is, the Hessian matrices $\nabla^2 f(x)$, $\nabla^2 g_i(x)$, and $\nabla^2 h_i(x)$. For nonlinear optimization, AMLs take closed-form algebraic equations and automatically generate routines for exact derivative evaluations in a form that solvers may call directly. An alternative to using AMLs here is to use general-purpose automatic differentiation (AD) tools which can be used to evaluate derivatives of code, an option which will be further discussed in Section~\ref{section:nonlinear}.
Lacking AMLs or other AD tools, one is faced with the tedious and error-prone task of implementing code to evaluate derivatives manually~\cite[p. 297]{gill1981practical}. In such cases it is common to forgo second-order derivatives or even first-order derivatives, even when providing them could reduce the solution time.

While the computational efficiency of an AML's translation of user input into solver input is something that can be empirically measured, we must note that the intrinsic usefulness of an AML is derived from how ``naturally'' the original mathematical statement can be translated into code.
This is a more subjective proposition, so in Figure~\ref{fig:codeexamples} we present the formulation in JuMP, AMPL, Pyomo and GAMS of a minimum cost flow problem as a linear program (see, e.g., \cite{bertsimas1997introduction}) over a graph $G=(V,E)$, where the vertices $V = \{1,2,\ldots,n\}$ are consecutively numbered with a ``source'' at vertex 1 and ``sink'' at vertex $n$:
\begin{equation}\label{mincostflow}
\begin{alignedat}{2}
\color{objcol}{\blacksquare} & \quad &
\min_{x}\quad & \sum_{(i,j)\in E} c_{i,j} x_{i,j} \\
\color{flowconcol}{\blacksquare} & \quad &
\text{s.t.}\quad
& \sum_{(i,j)\in E} x_{i,j}
    = \sum_{(j,k)\in E} x_{j,k} \quad j = 2, \ldots, n-1\\
\color{totflowcol}{\blacksquare} & \quad &
& \sum_{(i,n)\in E} x_{i,n} = 1 \\
\color{vardefcol}{\blacksquare} & \quad &
& 0 \leq x_{i,j} \leq C_{i,j} \quad \forall (i,j)\in E
\end{alignedat}
\end{equation}
The four AMLs share much in common: all involve declaring a set of variables indexed by the set of edges, all have a line for setting the objective function, and all have methods for iterating over a range ($2$ to $n-1$) and taking a sum over variables subject to some condition.
Of the four, JuMP and AMPL are perhaps the most similar, although JuMP benefits from being embedded in a full programming language, allowing us to define an Edge type that stores the problem data in a succinct and intuitive fashion. Pyomo is also embedded in a programming language, but as Python doesn't have the same syntactic macro functionality as Julia, some things are more uncomfortable than is ideal (setting the variable upper bounds, indexed constraint construction). Finally GAMS has perhaps the most verbose and idiosyncratic syntax, with features like set filtering with the \$ character that are not commonly found in either programming or modeling languages.
Our main claim is that JuMP is a ``natural'' and easy-to-use modeling language, and for the rest of this paper will instead focus on the technical details that allow it to be efficient and to enable unique features not found in other AMLs.

\begin{figure}
    \centering
    
    \textbf{JuMP}

\begin{Verbatim}[commandchars=\\\{\},codes={\catcode`$=3\catcode`^=7}]
{\color{white}$\blacksquare$} immutable Edge
{\color{white}$\blacksquare$}     from; to; cost; capacity
{\color{white}$\blacksquare$} end
{\color{white}$\blacksquare$} edges = [Edge(1,2,1,0.5), Edge(1,3,2,0.4), Edge(1,4,3,0.6),
{\color{white}$\blacksquare$}          Edge(2,5,2,0.3), Edge(3,5,2,0.6), Edge(4,5,2,0.5)]
{\color{white}$\blacksquare$} mcf = Model()
{\color{vardefcol}$\blacksquare$} @variable(mcf, 0 <= flow[e in edges] <= e.capacity)
{\color{totflowcol}$\blacksquare$} @constraint(mcf, sum\{flow[e], e in edges; e.to==5\} == 1)
{\color{flowconcol}$\blacksquare$} @constraint(mcf, flowcon[n=2:4], sum\{flow[e], e in edges; e.to==node\}
{\color{flowconcol}$\blacksquare$}                                  == sum\{flow[e], e in edges; e.from==node\}) 
{\color{objcol}$\blacksquare$} @objective(mcf, Min, sum\{e.cost * flow[e], e in edges\})
\end{Verbatim}

    \vspace{0.2cm}
    \textbf{AMPL}

\begin{Verbatim}[commandchars=\\\{\},codes={\catcode`$=3\catcode`^=7}]
{\color{white}$\blacksquare$} set edges := \{(1,2),(1,3),(1,4),(2,5),(3,5),(4,5)\};
{\color{white}$\blacksquare$} param cost\{edges\};  param capacity\{edges\};
{\color{white}$\blacksquare$} data ...; # Data is typically stored separately in AMPL;
{\color{vardefcol}$\blacksquare$} var flow\{(i,j) in edges\} >= 0.0, <= capacity[i,j];
{\color{totflowcol}$\blacksquare$} subject to unitflow:  sum\{(i,5) in edges\} flow[i,5] == 1;
{\color{flowconcol}$\blacksquare$} subject to flowconserve \{n in 2..4\}:
{\color{flowconcol}$\blacksquare$}   sum\{(i,n) in edges\} flow[i,n] == sum\{(n,j) in edges\} flow[n,j];
{\color{objcol}$\blacksquare$} minimize flowcost:    sum\{(i,j) in edges\} cost[i,j] * flow[i,j];
\end{Verbatim}
    
    \vspace{0.2cm}
    \textbf{Pyomo}
    
\begin{Verbatim}[commandchars=\\\{\},codes={\catcode`$=3\catcode`^=7}]
{\color{white}$\blacksquare$} edges    = [(1,2),     (1,3),     (1,4),     (2,5),     (3,5),     (4,5)]
{\color{white}$\blacksquare$} cost     = \{(1,2):1,   (1,3):2,   (1,4):3,   (2,5):2,   (3,5):2,   (4,5):2\}
{\color{white}$\blacksquare$} capacity = \{(1,2):0.5, (1,3):0.4, (1,4):0.6, (2,5):0.3, (3,5):0.6, (4,5):0.5\}
{\color{white}$\blacksquare$} mcf = ConcreteModel()
{\color{vardefcol}$\blacksquare$} mcf.flow = Var(edges, bounds=lambda m,i,j: (0,capacity[(i,j)]))
{\color{totflowcol}$\blacksquare$} mcf.uf = Constraint(expr=sum(mcf.flow[e] for e in edges if e[1]==5) == 1)
{\color{flowconcol}$\blacksquare$} def con_rule(mcf,n): return sum(mcf.flow[e] for e in edges if e[1]==n) ==
{\color{flowconcol}$\blacksquare$}                             sum(mcf.flow[e] for e in edges if e[0]==n)
{\color{flowconcol}$\blacksquare$} mcf.flowcon = Constraint([2,3,4],rule=con_rule)
{\color{objcol}$\blacksquare$} mcf.flowcost = Objective(expr=sum(cost[e]*mcf.flow[e] for e in edges))
\end{Verbatim}
    
    \vspace{0.2cm}
    \textbf{GAMS}
    
\begin{Verbatim}[commandchars=\\\{\},codes={\catcode`$=3\catcode`^=7}]
{\color{white}$\blacksquare$} SET nodes /n1*n5/; SET midnodes(nodes) /n2*n4/; SET lastnode(nodes) /n5/;
{\color{white}$\blacksquare$} ALIAS(nodes,nodefrom,nodeto,n);
{\color{white}$\blacksquare$} SET edges(nodes,nodes) / n1.n2 n1.n3 n1.n4 n2.n5 n3.n5 n4.n5 /;
{\color{white}$\blacksquare$} PARAMETER cost(nodes,nodes) / ... /;      * Data omitted
{\color{white}$\blacksquare$} PARAMETER capacity(nodes,nodes) / ... /;  * for space reasons
{\color{vardefcol}$\blacksquare$} POSITIVE VARIABLE flow(nodefrom,nodeto); flow.UP(edges) = capacity(edges);
{\color{totflowcol}$\blacksquare$} EQUATION unitflow;
{\color{totflowcol}$\blacksquare$} unitflow.. sum\{edges(nodefrom,lastnode), flow(nodefrom,lastnode)\} =e= 1;
{\color{flowconcol}$\blacksquare$} EQUATION flowcon(nodes);
{\color{flowconcol}$\blacksquare$} flowcon(midnodes(n)).. sum\{edges(nodefrom,n), flow(nodefrom,n)\} =e=
{\color{flowconcol}$\blacksquare$}                sum\{edges(n,nodeto), flow(n,nodeto)\};
{\color{objcol}$\blacksquare$} FREE VARIABLE obj;
{\color{objcol}$\blacksquare$} EQUATION flowcost; flowcost.. obj =e= sum\{edges, cost(edges)*flow(edges)\};
{\color{white}$\blacksquare$} MODEL mincostflow /all/; SOLVE mincostflow USING lp MINIMIZING obj;
\end{Verbatim}

    \caption{Modeling a minimum cost flow problem in JuMP, AMPL, Pyomo, and GAMS. The colored squares show the correspondence between the code and the four components of Equation~\ref{mincostflow}. For concreteness, we provide an explicit example of a five-node problem with data when it fits. The JuMP and Pyomo examples are complete, valid code (as of this writing) and can be copy-pasted into a terminal to run after importing the corresponding packages.}
    \label{fig:codeexamples}
\end{figure}

The technical tasks that an AML must perform can be roughly divided into two simple categories: first, to load the user's input into memory, and second, to generate the input required by the solver, according to the class of the problem. For both of these tasks, we have made some unorthodox design decisions in JuMP in order to achieve good performance under the constraints of being embedded within a high-level language. We will review these in the following sections.

We note that, as we have mentioned, JuMP provides access to a number of advanced techniques which have not been typically available in AMLs. For example, \textit{branch-and-cut} is a powerful technique in integer programming for accelerating the solution process by dynamically improving the convex (linear) relaxations used within the branch-and-bound algorithm. Users wanting to extend a solver's branch-and-cut algorithm with dynamically generated ``cuts'' for a particular problem structure have typically needed to turn to low-level coding in C++ for an efficient implementation via callback functions, since this interaction requires bidirectional communication with a solver during the solution process. To our knowledge, JuMP is the first AML to provide a simple, high-level, and efficient (in-memory) interface to branch-and-cut and other similar techniques.
This feature has already seen fruitful use in research~\cite{nikita,VielmaExtendedFormulations} and teaching~\cite{IAPClass}.

\section{Syntactic macros: parsing without a parser}\label{sec:macros}

AMLs like AMPL and GAMS are standalone in the sense that they have defined their own syntax entirely separate from any existing programming language. They have their own formats for providing input data (although they can also connect to databases and spreadsheets) and implement custom parsers for their proprietary syntax; for example, AMPL uses the LEX and YACC parser generator utilities~\cite{amplpaper}.

Embedding an AML within an existing programming language brings with it the benefit of being able to bootstrap off the existing, well defined grammar and syntax of the language, eliminating a complex part of implementing an AML. Perhaps more importantly for users, embedded AMLs typically allow interlacing the AML's math-like statements declaring an optimization problem with arbitrary code which may be used to prepare input data or process the results. However, embedding also brings with it the challenge of obtaining the desired expressiveness and ease of use within the limits of the syntax of the parent language.

The most common approach (taken by Pyomo, YALMIP, and others) to capturing user input is \textit{operator overloading}. One introduces a new class of objects, say, to represent a decision variable or vector of decision variables, and extends the language's definition of basic operators like $+$, $*$, $-$, etc, which, instead of performing arithmetic operations, build up data structures which represent the expression. For example, to represent a quadratic expression $\sum_{(i,j) \in J} b_{ij} x_i x_j + \sum_{i \in I} a_i x_i + c$, one stores the constant $c$, the coefficient vectors $b$, $a$, and the index sets $I$ and $J$. Letting $n$ be the number of decision variables in a problem, an unfortunate property of addition of two quadratic expressions is that the size of the resulting expression is not bounded by a constant independent of $n$, simply because the coefficient and index vectors can have as many as $O(n^2)$ terms. This means that basic operations like addition and subtraction are no longer fast, constant-time operations, a property which is almost always taken for granted in the case of floating-point numbers.  As a concrete example, consider the following quadratic expression in the variable $x$ indexed over $\{1,\ldots,d\} \times \{1,\ldots,d\}$:
\begin{equation}\label{eq:quadexample}
    1 + \sum_{i=1}^d\sum_{j=1}^d |c_j - i|(1-x_{i,j})x_{1,j}
\end{equation}

In Python, one might naturally express~\eqref{eq:quadexample} as
\begin{lstlisting}
1 + sum(abs(c[j]-i)*(1-x[i,j])*x[0,j] for i in range(d) for j in range(d))
\end{lstlisting}

which takes advantage of the built-in \texttt{sum} command which internally accumulates the terms one-by-one by calling the addition operator $d^2$ times. The partial sums generated with each addition operation are quadratic expressions which have $O(d^2)$ terms, so this naive approach can have a cost of $O(d^4) = O(n^2)$ operations and excessive memory allocations. An obvious workaround for this issue is to accumulate the terms in a single output expression instead of a generating a new expression for each partial sum. While there are a number of ways to mitigate this slow behavior within the framework of operator overloading, our benchmarks will demonstrate that they may not be sufficient to achieve the best performance.

When designing JuMP, we were not satisfied by the performance limitations of operator overloading and instead turned to an advanced feature of Julia called \textit{syntactic macros}~\cite{JuliaArxiv}.
Readers may be familiar with macros in C and C++ which perform textual substitutions; macros in Julia are much more powerful in that they function at the level of syntax. For example, the expression~\eqref{eq:quadexample} could be written in JuMP syntax as
\begin{lstlisting}
@expression(1 + sum{abs(c[j]-i)*(1-x[i,j])*x[1,j], i in 1:N, j in 1:N})
\end{lstlisting}

The \texttt{@} sign denotes a call to a macro named \texttt{expression}, which constructs a JuMP expression object. The input to the macro will be a data structure representing the \textit{Julia expression} contained within, not simply a string of text. That is, Julia's internal parser will be invoked to parse the expression, but instead of directly evaluating it or compiling it to code, it will be sent to a routine written in Julia which we (as authors of JuMP) have defined. Note that the syntax \texttt{sum\{\}} is generally not valid Julia code for computing a sum, although it is recognized by the Julia parser because the syntax is used in other contexts, which allows us to endow this syntax with a new meaning in the context of JuMP\footnote{As of Julia 0.5, JuMP will transition to the new syntax \texttt{sum(abs(c[j]-i)*(1-x[i,j])*x[1,j] for i in 1:N, j in 1:N)} which will more closely match Julia code.}.

Macros enable JuMP to provide a natural syntax for algebraic modeling without writing a custom text-based parser and without the drawbacks of operator overloading. Within the computer science community, macros have been recognized as a useful tool for developing domain-specific languages, of which JuMP is an example~\cite{dsl}. Indeed, the implementation of macros in Julia draws its inspiration from Lisp~\cite{JuliaArxiv}. However, such functionality historically has not been available within programming languages targeted at scientific computing, and, to our knowledge, JuMP is the first AML to be designed around syntactic macros.

\section{Code generation for linear and conic-quadratic models}\label{sec:linear}
Linear and conic-quadratic optimization problems are essential and surprisingly general modeling paradigms that appear throughout operations research and other varied fields---often at extremely large scales. Quadratic optimization generalizes linear optimization by allowing convex quadratic terms $\frac{1}{2}x^TQx$ in the objective and constraints. Conic-quadratic, also known as second-order cone, optimization generalizes quadratic optimization with constraints of the form $||x||_2 \leq t$,
where both $x$ and $t$ are decision variables~\cite{socp}.
Computational tools for \emph{solving} these problems derive their success from exploiting the well-defined structure of these problems. Analogously, JuMP is able to efficiently process large-scale problems by taking advantage of structural properties and generating efficient code through Julia's code generation functionality.

Julia is, at the same time, both a dynamic and compiled language. Julia uses the LLVM compiler~\cite{lattner2004llvm} dynamically at runtime, and can generate efficient, low-level code as needed. This technical feature is one of the reasons why Julia can achieve C-like performance in general~\cite{Julia}, but we will restrict our discussion to how JuMP takes advantage of it.

In the previous section we described how JuMP uses macros to accept user input in the form of a data structure which represents the input expression. The other side of macros is \textit{code generation}. More specifically, macros can be understood as functions whose input is code and whose output is code. Given an input expression, a macro produces a data structure which represents an output expression, and that expression is then substituted in place and compiled. For example, the call to the \texttt{expression} macro in Section~\ref{sec:macros} would output, in pseudo-code from, the following code:

\begin{lstlisting}
Initialize an empty quadratic expression q
Add 1 to the constant term
Count the number of terms K in the sum{} expression
Pre-allocate the coefficient and index vectors of q to hold K elements
for i in 1:d, j in 1:d
    Append -abs(c[j]-i)*x[i,j]*x[1,j] to the quadratic terms in q
    Append abs(c[j]-i)*x[1,j] to the linear terms in q
end
\end{lstlisting}

Note that this code runs in $O(d^2)$ operations, a significant improvement over the $O(d^4)$ naive operator overloading approach. The code produced is also similar to a hand-written matrix generator. Indeed, one could summarize JuMP's approach to generating linear and quadratic models as translating users' algebraic input into fast, compiled code which acts as a matrix generator. JuMP's approach to semidefinite optimization, a recently added feature which we will not discuss further, generally follows this path but also employs operator overloading for certain matrix operations.

\subsection{Benchmarks}\label{sec:conicbenchmark}

We now provide computational evidence that JuMP is able to produce quadratic and conic-quadratic optimization models, in a format suitable for consumption by a solver, as fast as state-of-the-art commercial modeling languages.
To do so we measure the time elapsed between launching the executable that builds the model and the time that the solver begins the solution process, as determined by recording when the first output appears from the solver. 
This methodology allows the modeling language to use a direct in-memory solver interface if it desires, or in the case of some tools a compact file representation.
We selected Gurobi~6.5.0~\cite{gurobi} as the solver, and evaluated the following modeling systems: the Gurobi C++ interface (based on operator overloading), JuMP~0.12 with Julia~0.4.3, AMPL~20160207 \cite{AMPLBook}, GAMS~24.6.1 \cite{GAMS}, Pyomo~4.2.10784 with Python~2.7.11 \cite{Pyomo}, and CVX~2.1 \cite{cvx-matlab} and YALMIP~20150918 \cite{YALMIP} with MATLAB~R2015b. These particular modeling systems are chosen for being widely used in practice within different communities. The benchmarks were run on a Linux system with an Intel Xeon CPU E5-2650 processor.

We implemented two different optimization problems in all seven modeling languages: a linear-quadratic control problem (\texttt{lqcp}) and a facility location problem (\texttt{fac}). We do not claim that these models are representative of all conic-quadratic problems; nevertheless, they provide a good stress test for generating models with many quadratic and conic-quadratic terms. The models are further described in the appendix. Models using the C++ interface are implemented in a way that mitigates the drawbacks of operator overloading by appending to existing expressions using the \texttt{+=} operator; such approaches, however, are not idiomatic in Pyomo, CVX, or YALMIP. The results (Table~\ref{tab:linearbench}) show that for \texttt{lqcp}, JuMP, AMPL, and the C++ interface are roughly equivalent at the largest scale, with GAMS and YALMIP approximately four times slower and CVX thirteen times slower than JuMP. Pyomo was significantly slower and was unable to construct the largest model within ten minutes. For \texttt{fac}, JuMP, AMPL, GAMS and the C++ interface times all perform roughly the same, while Pyomo is unable to build the largest instance with ten minutes, YALMIP can build only the smallest instance within the time limit, and CVX is unable to build any instances within the time limit. These results demonstrate that JuMP can be reasonably competitive with widely used commercial systems and in some cases significantly faster than open-source alternatives.

\begin{table}
\centering
\begin{tabular}{c | r | r r r | r r r}
 \multicolumn{2}{c}{}  & \multicolumn{3}{c}{\textbf{Commercial}} & \multicolumn{3}{c}{\textbf{Open-source}} \\
Instance  & \textbf{JuMP} & GRB/{\footnotesize C++} & AMPL & GAMS & Pyomo & CVX & YALMIP \\ \hline
lqcp-500  &  8 &  2 &  2 &   2 &               55 &  6 &   8 \\
lqcp-1000 & 11 &  6 &  6 &  13 &              232 &  48 &  25 \\
lqcp-1500 & 15 & 14 & 13 &  41 &              530 &  135  &  52 \\
lqcp-2000 & 22 & 26 & 24 & 101 & \textgreater 600 & 296  & 100 \\[3pt]
fac-25  &  7 &  0 &  0 &  0 &               14 & \textgreater 600  &              533 \\
fac-50  &  9 &  2 &  2 &  3 &              114 & \textgreater 600  & \textgreater 600 \\
fac-75  & 13 &  5 &  7 & 11 &              391 & \textgreater 600  & \textgreater 600 \\
fac-100 & 24 & 12 & 18 & 29 & \textgreater 600 & \textgreater 600  & \textgreater 600 \\
\end{tabular}

\caption{Time (sec.) to generate each model and pass it to the solver, a comparison between JuMP and existing commercial and open-source modeling languages. The \texttt{lqcp} instances have quadratic objectives and linear constraints. The \texttt{fac} instances have linear objectives and conic-quadratic constraints.}
\label{tab:linearbench}
\end{table}

\subsection{Optimizing a sequence of models}
As we observe in Table \ref{tab:linearbench}, JuMP has a noticeable start-up cost of a few seconds even for the smallest instances. This start-up cost is primarily composed of compilation time; however, if a family of models is solved multiple times within a single session, this cost of  compilation is only paid for the \textit{first} time that an instance is solved. That is, when solving a sequence of instances in a loop, the amortized cost of compilation is negligible.

Indeed, solving a sequence of related optimization problems is a common idiom when performing exploratory analysis or implementing more advanced algorithms. A number of algorithms including branch-and-bound~\cite{Atamturk:2004}, Benders decomposition~\cite{StochasticBenders}, and cutting-plane methods derive their efficiently from the fact that when solving a sequence of \textit{linear} programming problems, one can ``hot-start'' the simplex algorithm from the previous optimal solution when new constraints are added or when objective or right-hand side coefficients are modified. JuMP supports all of these classes of modifications by enabling
efficient in-memory access to solvers so that they can maintain their internal state across solves, when possible, avoiding the significant overhead that would be incurred when regenerating a model from scratch inside of a loop (even if information on the previous optimal solution persists). For example, a straightforward implementation of a cutting-plane algorithm in GAMS was found to be 5.8x slower overall than an implementation of the same algorithm in C++~\cite{GAMSGUSS}, illustrating the overhead induced by these traditional AMLs which do not allow direct in-memory access to solvers.

Here we present a small example that demonstrates solving and modifying a problem in a loop with JuMP,

\begin{lstlisting}
l2approx = Model()
@variable(l2approx, -1 <= x[1:N] <= +1)
@objective(l2approx, Max, sum(x))
solve(l2approx); v = getvalue(x)  # Build and solve for initial solution
while norm(v) >= 1 + tol
    @constraint(l2approx, dot(v,x) <= norm(v))
    solve(l2approx); v = getvalue(x)  # Optimize from prev solution
end
\end{lstlisting}
We maximize a simple linear function ($\sum_{i=1}^N x_i$) over the $\ell_2$ ``ball'' constraint $\Vert x \Vert_2 \leq 1$ by approximating this nonlinear constraint with a finite sequence of tangent hyperplanes generated only as needed, allowing us to use an LP solver instead of a solver that supports conic-quadratic constraints. Direct extensions of this technique have proven useful for solving mixed-integer conic optimization problems~\cite{VielmaExtendedFormulations}.
When a constraint is added to the JuMP model by the user, as in the above example, JuMP adds the constraint directly to the solver's in-memory representation of the problem, rather than generating the whole model from scratch. As a result the solver is able to use the previous (now infeasible) solution as a hot-start by applying the dual simplex method.

JuMP's design comes for the price of not supporting constructs for parametric data~\cite{AMPLBook} as GAMS, AMPL, and Pyomo do; that is, in JuMP one cannot define a parametric value and have its values propagate automatically through an LP model as the value of the parameter changes, because doing so would complicate the abstraction layer between JuMP's representation of a model and the solver's in-memory representation of a model\footnote{In the case of nonlinear optimization, JuMP offers parameter objects because these can be efficiently integrated within the derivative computations.}.
For more complex changes to problem data and structure, such as modifying the coefficients in existing constraints, the idiomatic approach when using JuMP is to construct a new model from scratch, possibly inside a function that takes parameters as its input. The following pseudocode, with JuMP on the left and AMPL on the right, demonstrates the difference in approach from a user's perspective.

\begin{minipage}{.5\textwidth}
\begin{lstlisting}
function f(p)
    # Build and solve model
    # using p as data, then
    # process solution
end
f(1)
f(2)
f(3)
\end{lstlisting}
\end{minipage} 
\begin{minipage}{.5\textwidth}
\begin{lstlisting}
param p;
# Define algebraic model using p
let p := 1;
solve;  # ... process solution
let p := 2;
solve;  # ... process solution
let p := 3;
solve;  # ... process solution
\end{lstlisting}
\end{minipage}

\section{Computing derivatives for nonlinear models}
\label{section:nonlinear}
Recall that the role of a modeling language for nonlinear optimization is to allow users to specify ``closed-form,'' algebraic expressions for the objective function $f$ and constraints $g$ and $h$ in the formulation~\eqref{eq:nlp-standard} and communicate first-order and typically second-order derivatives with the optimization solver. Commercial modeling languages like AMPL and GAMS represent the state of the art in modeling languages for nonlinear optimization. Likely because of the increased complexity of computing derivatives, even fewer open-source implementations exist than for linear or quadratic models.

Notable alternative approaches to traditional algebraic modeling for nonlinear optimization include CasADi~\cite{CasADi} and CVX~\cite{cvx-matlab}. CasADi allows interactive, scalar- and matrix-based construction of nonlinear expressions via operator overloading with automatic computation of derivatives for optimization. CasADi has specialized features for optimal control but, unlike traditional AMLs, does not support linear optimization as a special case. CVX, based on the principle of disciplined convex programming (DCP)~\cite{DCP}, allows users to express convex optimization problems in a specialized format which can be transformed into or approximated by conic programming without the need for computing derivatives\footnote{Note that the DCP paradigm is available in Julia through the Convex.jl package~\cite{convexjl}.}. The traditional nonlinear optimization formulation considered here applies more generally to derivative-based convex and nonconvex optimization. 

JuMP, like AMPL and GAMS, uses techniques from automatic (or algorithmic) differentiation (AD) to evaluate derivatives of user-defined expressions. 
In this section, we introduce these techniques with a focus on how AMLs relate to more general-purpose AD tools. In this vein, we discuss JuMP's unique ability to automatically differentiate \textit{user-defined functions}. Concluding this section, we present a set of performance benchmarks.

\subsection{Expression graphs and reverse-mode AD}
\begin{figure}[t]
        \centering
        \raisebox{-\height}{
\begin{tikzpicture}[scale=1.0]
        \tikzset{every internal node/.style={circle,draw}}
        \tikzset{every tree node/.style={font=\tt}}
        \tikzset{every leaf node/.style={rectangle,draw,minimum height=18pt}}
        \Tree [.$\exp(\cdot)$ [.+ [.$(\cdot)^2$ [.$x$ ] ] [.$(\cdot)^2$ [.$y$ ] ] ] ]
\end{tikzpicture}
}
\caption{A graph representation of the nonlinear expression $\exp(x^2+y^2)$. JuMP uses this expression graph structure for efficient evaluation of derivatives.}
\label{fig:expgraph}

\end{figure}
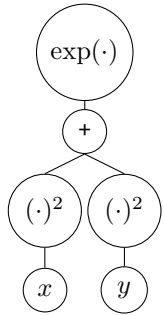

    A natural data structure for representing nonlinear expressions is the \textit{expression graph}, which is a directed acyclic graph (or typically, a tree) that encodes the sequence of operations required to compute the expression as well as the dependency structure between operations. For example, Figure~\ref{fig:expgraph} illustrates how the nonlinear expression $\exp(x^2+y^2)$ is represented as a graph, with nodes representing the input values $x$ and $y$ together with every ``basic'' operation like addition and exponentiation that is performed in computing the value of the expression. We will return later to the question of what is considered a basic operation, but for now consider these as all operations that one might compose to write down a closed-form algebraic equation. Edges in the expression graph represent immediate dependencies between operations. The expression graph encodes all needed information for JuMP to evaluate and compute derivatives of nonlinear expressions, and JuMP generates these objects by using macros analogously to how JuMP generates sparse matrix data structures representing linear and quadratic functions.

    While general-purpose AD tools like ADOL-C~\cite{ADOLC} use operator overloading (or direct analysis of source code~\cite{ADIFOR}) to generate these expression graphs from arbitrary code; AMLs like AMPL, GAMS, and JuMP have an easier task because user input is constrained to follow a specific syntax and thus these AMLs are generally more reliable. The value of using JuMP, for example, over using more general-purpose AD tools is that JuMP provides a guarantee to the user that all input following its simple syntax can be differentiated efficiently, the only limitation being that the expression graph objects must fit within memory. On the other hand, making use of more general-purpose tools requires a nontrivial amount of expertise (for example, preparing C++ code for ADOL-C requires extensive modifications and use of specialized assignment operators). We have found users from fields like statistics who have traditionally not been users of AMLs being drawn to JuMP for its AD features and its being embedded in a familiar programming language, to the extent that they are willing to rewrite complex statistical models into JuMP's syntax~\cite{Giordano}. It remains to be seen how general this adoption trend may be, but we believe that there is large scope for judicious use of AMLs as AD tools within domains that have not widely adopted them so far.

    Given an expression graph object, one can compute the numerical value of the expression by iterating through the nodes of the graph in an order such that by the time we reach a given node to evaluate its corresponding operation, the numerical values of all its inputs (children) have already been computed. A perhaps surprising result is that it is possible to apply the chain rule in such a way that by iterating through the nodes in the reverse order (parents before children), in a single pass, we obtain the exact gradient vector $\nabla f(x)$. This reverse-pass algorithm, known suitably as \textit{reverse-mode AD}, delivers gradients of $f$ for a small constant factor times the cost evaluating $f$ itself. We refer readers to~\cite{Griewank2008EDP,NaumannAD} for further discussion of this powerful technique which JuMP employs.

\subsection{User-defined functions}
    Modeling languages like JuMP include a basic library of nonlinear functions which are available for use within expressions. JuMP's basic library is extensive and includes special functions like the error function $\operatorname{erf}()$, which enables users to express Gaussian cumulative densities, for example. AMPL recently developed an interface to functions from the GNU Scientific Library~\cite{GSL}, greatly extending the range of available functions. However, in the cases where the built-in library is insufficient, historically there has been no user-friendly way to incorporate \textit{user-defined functions} into AMLs.\footnote{Doing so within AMPL, for example, requires developing a shared library object in C which links to the AMPL solver library.} A compelling application for user-defined functions is optimal control problems constrained by differential equations, where standalone integrators are used to enforce dynamic constraints~\cite{CasADiDynamic}. JuMP is the first AML to provide a simple interface not only for user-defined functions with user-provided (hand-coded) derivatives but also to provide an option to \textit{automatically differentiate user-defined functions}. We provide a brief example below of this usage in JuMP and then describe the implementation. 

\begin{lstlisting}
function squareroot(x)
    z = x # Initial starting point for Newton's method
    while abs(z*z - x) > 1e-13
        z = z - (z*z-x)/(2z)
    end
    return z
end
JuMP.register(:squareroot, 1, squareroot, autodiff=true)

m = Model()
@variable(m, x[1:2], start=0.5)
@objective(m, Max, sum(x))
@NLconstraint(m, squareroot(x[1]^2+x[2]^2) <= 1)
solve(m)
\end{lstlisting}

First, we define the \texttt{squareroot} function using generic Julia code. This function computes the square root of a number by applying Newton's method to find the zero of the function $f(z) = z^2 - x$. The function \texttt{JuMP.register} registers the nonlinear function with the symbolic name \texttt{squareroot} and passes a reference to the function defined above. The second argument \texttt{1} indicates that the input to the function is univariate. The \texttt{autodiff=true} option instructs JuMP to automatically compute the derivatives of the user-defined function (if this option is not set, users must also provide derivative evaluation callbacks). The subsequent JuMP model defines the problem of maximizing $x_1 + x_2$ subject to the constraint $\sqrt{x_1^2 + x_2^2} \le 1$, where the square root is computed by the user-defined Julia code. 

In principle, JuMP could apply reverse-mode AD to user-defined functions by using operator overloading to build up the expression graph representation of the function, using analogous techniques to ADOL-C. However, no mature implementation of this approach currently exists in Julia. Instead, JuMP uses ForwardDiff.jl, a standalone Julia package implementing \textit{forward-mode} AD~\cite{forwarddiff}, to compute derivatives of these functions.

Forward-mode AD, known as such because it can be applied with a single forward pass though the expression graph data structure, can be interpreted as computing directional derivatives by introducing infinitesimal perturbations~\cite{Griewank2008EDP}. In comparison, the more well-known finite differencing method employs small but finite perturbations. The operator overloading approach for forward-mode AD is
to introduce a new class of number $a+b\epsilon$ where $\epsilon^2 = 0$ (analogously to $i^2=-1$ for the complex numbers). 
The implementation in Julia's ForwardDiff.jl is conceptually quite similar to that in other languages, and we refer readers to Neidinger~\cite{IntroAD} for a comprehensive introduction to forward-mode AD and its implementation in MATLAB using operator overloading. In Julia, however, user-defined numerical types are given first-class treatment by the compiler and produce efficient low-level machine code, which is not the case for MATLAB. Note that forward-mode applied in this way does not require an explicit expression graph representation of a function and hence is simpler to implement than reverse-mode.

The only burden on users providing functions for automatic differentiation is to write code which is generic with respect to the type of the numerical input (in the above example, \texttt{x}). This design is equivalent in spirit to using templates in C++ but with a much less heavyweight syntax. The cost of applying forward-mode AD grows linearly with the the input dimension of the target function; hence for high-dimensional user-defined functions, users may still want to provide derivative evaluation callbacks if speed is a concern. Nevertheless, the ability to automatically differentiate user-defined functions begins to blur the distinction between AMLs and more traditional AD tools, and we look forward to seeing the applications of this recently added feature.

\subsection{From gradients to Hessians}
In addition to gradients, off-the-shelf nonlinear optimizers typically request second-order derivatives. A basic operation for computing second-order derivatives is the Hessian-vector product $\nabla^2 f(x)d$. Since this product is a directional derivative of the gradient, we now have the tools to compute it, by applying forward-mode automatic differentiation to the reverse-pass algorithm for gradients.

This composition of forward-mode with reverse-mode AD is known as forward-over-reverse mode~\cite{Griewank2008EDP,NaumannAD}, and JuMP implements it by manually augmenting the reverse-mode implementation to propagate the required infinitesimal components for directional derivatives. Note that we do not yet support second order automatic differentiation of user-defined functions.

Given this forward-over-reverse routine to compute Hessian-vector products, one could recover a dense Hessian matrix $\nabla^2 f(x)$ by calling the routine $n$ times, taking the $n$ distinct unit vectors. However, for large $n$, this method quickly becomes prohibitively expensive. By exploiting the sparsity structure of $\nabla^2 f(x)$, one instead may compute the entries of the Hessian matrix with far fewer than $n$ Hessian-vector products. For example, if the Hessian is known to be diagonal, one needs only a single Hessian-vector product with $d = (1,1,\cdots,1)^T$ to compute all nonzero elements of the Hessian. In general, the problem of choosing a minimal number of Hessian-vector products to compute all nonzero elements is NP-hard; we implement the acyclic graph coloring heuristic of Gebremedhin et al.~\cite{GebIjoc}. See Figure~\ref{fig:graphcolor} for an illustration. The Hessian matrices of typical nonlinear models exhibit significant sparsity, and in practice a very small number of Hessian-vector products are needed even for high-dimensional problems. We note that AMPL exploits Hessian structure through partial separability~\cite{AMPLad} instead of using graph coloring techniques.

\begin{figure}[t]
\centering
\[
\left[\begin{array}{>{\columncolor{red1!80}}c>{\columncolor{red2!80}}c>{\columncolor{red1!80}}c>{\columncolor{red2!80}}c>{\columncolor{red2!80}}c} h_{11} & h_{12} && h_{14}&\\
h_{12} & h_{22} & h_{23}&&\\
&h_{23} & h_{33}&& \\
h_{14} &&&h_{44}&\\
&&&&h_{55}
\end{array}\right]
\left[
\begin{array}{>{\columncolor{red1!80}}c>{\columncolor{red2!80}}c}
1 & \\
&1\\
1&\\
&1\\
&1
\end{array}\right]
 =
\left[
\begin{array}{cc}
h_{11} & h_{12} + h_{14}\\
h_{12} +h_{23} & h_{22}\\
h_{33} & h_{23}\\
h_{14} & h_{44}\\
& h_{55}
\end{array}\right]
\]
\caption{
Many solvers can benefit from being provided the Hessian matrix of second-order derivatives at any point. JuMP uses reverse-mode automatic differentiation to generate a ``black box'' routine that computes \textit{Hessian-vector products} and uses this to calculate the non-zero elements of the Hessian matrix. For efficiency, we would like to use as few Hessian-vector products as possible; by using a specialized graph coloring heuristic \cite{GebIjoc}, we can find a small number of evaluations to do so. Above, we illustrate a symmetric $5 \times 5$ Hessian matrix with $h_{ij} = \frac{\partial^2 f}{\partial x_i\partial x_j}(x)$ for some $f$. Omitted entries are known to be zero. In this example, only two Hessian-vector products are needed.
}
\label{fig:graphcolor}
\end{figure}

\subsection{Benchmarks}\label{sec:nlbenchmark}

We now present benchmarks evaluating the performance of JuMP for modeling nonlinear optimization problems. Similar to the experimental design in Section~\ref{sec:conicbenchmark}, we measure the time elapsed after starting the executable until the solver, Ipopt~\cite{ipopt}, reports the problem dimensions as confirmation that the instance is loaded in memory. Then, we fix the total number of iterations performed to three and record the time spent in function or derivative evaluations as reported by Ipopt. We evaluated the following modeling systems: JuMP, AMPL, Pyomo, GAMS, and YALMIP. Recall that CVX does not support derivative-based nonlinear models. Also, YALMIP does not support Hessian evaluations, so we measure only model generation time.

We test two families of problems, nonlinear beam control (\texttt{clnlbeam}) and nonlinear optimal power flow (\texttt{acpower}), which are further described in the appendix. These two problems stress different aspects of model generation; the \texttt{clnlbeam} family has large, sparse, and very simple nonlinear expressions with a diagonal Hessian matrix, while the \texttt{acpower} family has a smaller number of variables but much more complex nonlinear network structure with a Hessian matrix with an irregular sparsity pattern. For model generation times (Table~\ref{tab:nonlinearmodel}), JuMP has a relatively large startup cost, which is dominated by Julia's compiler. However, as the size of the instance increases, JuMP becomes significantly faster than Pyomo and YALMIP.
As suggested by its performance and the omission of Hessian computations, YALMIP's derivative-based nonlinear functionality is seemingly not designed for large-scale problems. We did not implement \texttt{acpower} in YALMIP.

The results in Table~\ref{tab:nonlineareval} compare the time spent evaluating derivatives. Excluding the smallest instances, JuMP remains within a factor of 2.2x of AMPL. JuMP is up to 3x faster than GAMS and in the worst case 25\% slower.
Note that Pyomo does not implement its own derivative computations; instead, it re-uses AMPL's open-source derivative evaluation library.

\begin{table}[ht]
\centering
\begin{tabular}{l|r|rr|rr}
\multicolumn{2}{c}{} & \multicolumn{2}{c}{\textbf{Commercial}} & \multicolumn{2}{c}{\textbf{Open-source}}\\
 Instance & \textbf{JuMP} & AMPL & GAMS & Pyomo & YALMIP\\ \hline
clnlbeam-5   & 12 &  0 &   0 &   5 & 76\\
clnlbeam-50  & 14 &  2 &   3 &  44 & \textgreater 600 \\
clnlbeam-500 & 38 & 22 &  35 & 453 & \textgreater 600 \\[3pt]
acpower-1    & 18 &  0 &   0 &   3 & - \\ 
acpower-10   & 21 &  1 &   2 &  26 & - \\ 
acpower-100  & 66 & 14 &  16 & 261 & -    
\end{tabular}
\caption{Time (sec.) to generate each model and pass it to the solver, a comparison between JuMP and existing commercial and open-source modeling languages for derivative-based nonlinear optimization. Dash indicates not implemented.}\label{tab:nonlinearmodel}
\end{table}

\begin{table}[ht]
\centering
\begin{tabular}{l|r|rr}
\multicolumn{2}{c}{} & \multicolumn{2}{c}{\textbf{Commercial}}\\
 Instance & \textbf{JuMP} & AMPL & GAMS\\\hline
clnlbeam-5   & 0.08 & 0.03 &  0.08 \\
clnlbeam-50  & 0.70 & 0.39 &  0.76 \\
clnlbeam-500 & 7.48 & 3.47 & 15.81 \\[3pt]
acpower-1    & 0.07 & 0.02 &  0.06 \\ 
acpower-10   & 0.66 & 0.30 &  0.53 \\ 
acpower-100  & 6.11 & 3.20 & 18.13    
\end{tabular}

\caption{Time (sec.) to evaluate derivatives (including gradients, Jacobians, and Hessians) during 3 iterations, as reported by Ipopt.  Pyomo relies on AMPL's ``solver library'' for derivative evaluations, and YALMIP does not provide second-order derivatives. }\label{tab:nonlineareval}
\end{table}

\section{Extensions}\label{sec:extensions}

JuMP is designed to be extensible, allowing for developers both to plug in new solvers for existing problem classes and to extend the syntax of JuMP itself to new classes of problems. In comparison, it is more common for AMLs to support only extending the set of solvers for existing, well-defined problem classes~\cite{SET}. A common thread motivating extensions to an AML's syntax, on the other hand, is that the more natural representation of a class of models may be at a higher level than a standard-form optimization problem. These classes of models furthermore may benefit from customized solution methods which are aware of the higher-level structure. Extensions to JuMP can expose these advanced problem classes and algorithmic techniques to users who just want to solve a model and not concern themselves with the low-level details.
We will present three extensions we recently developed with this motivation for handling different models for optimization under uncertainty:  parallel multistage stochastic programming, robust optimization, and chance constraints. While these three extensions were developed by the JuMP core developers, we would like to highlight that even more recently a number of syntactic extensions to JuMP have been developed independently~\cite{multijump,VI,complementarity}, illustrating the feasibility of doing so without intimate knowledge of JuMP's internals.

\subsection{Extension for parallel multistage stochastic programming}
The first example of a modeling extension built on top of JuMP is StructJuMP~\cite{StochJuMP} (formerly StochJuMP), a modeling layer for block-structured optimization problems of the form,
\begin{equation}\label{eq:2stageProb}
  \begin{array}{clcrcrcrcrcrcl}
\min
     & \multicolumn{12}{l}{\frac{1}{2}x_0^TQ_0x_0+c_0^Tx_0 +  \sum_{i=1}^N \left(\frac{1}{2}x_i^TQ_ix_i+c_i^Tx_i\right)}\\
\textnormal{s.t.}
     &  Ax_0     &           &           &            &         & = b_0,\\
     &  T_1x_0 + &  W_1x_1   &           &            &         & = b_1, \\
     &  T_2x_0 + &           & W_2x_2    &            &         & =  b_2, \\
     &\ \  \vdots&           &           &   \ddots   &         &\vdots\\
     &  T_Nx_0 + &           &           &            &  W_Nx_N & = b_N, \\
     &\multicolumn{12}{l}{x_0\geq 0,\quad x_1\geq 0,\quad x_2\geq 0,\quad\ldots,\quad x_N\geq 0.}
  \end{array}
\end{equation}

This structure has been well studied and arises from stochastic programming~\cite{Birge11}, contingency analysis~\cite{PhanPowerSystems}, multicommodity flow~\cite{PrimalBlockAngular}, and many other contexts. A number of specialized methods exist for solving problems with this structure (including the classical Benders decomposition method), and they require as input data structures the matrices $Q_i$, $T_i$, $W_i$, $A$, and vectors $c_i$ and $b_i$.

StructJuMP was motivated by the application to stochastic programming models for power systems control under uncertainty as outlined in \cite{Petra:2014}. For realistic models, the total number of variables may be in the tens to hundreds of millions, which necessitates the use of parallel computing to obtain solutions within reasonable time limits. In the context of high-performance computing, the phase of generating the model in serial can quickly become an execution bottleneck, in addition to the fact that the combined input data structures may be too large to fit in memory on a single machine. StructJuMP was designed to allow users to write JuMP models with annotations indicating the block structure such that the input matrices and vectors can be generated in parallel. That is, the entire model is not built in memory in any location: each computational node only builds the portion of the model in memory that it will work with during the course of the optimization procedure. This ability to generate the model in parallel (in a distributed-memory MPI-based~\cite{MPI} fashion) distinguishes StructJuMP from existing tools such as PySP~\cite{PySP}. 

StructJuMP successfully scaled up to 2048 cores of a high-performance cluster, and in all cases the overhead of model generation was a small fraction of the total solution time. Furthermore, StructJuMP was \textit{easy} to develop, consisting of less than 500 lines of code in total, which includes interfacing with a C++-based solver and the MPI library for parallel computing. For comparison, SML~\cite{sml}, an AMPL extension for conveying similar block structures to solvers, was implemented as a pre- and post-processor for AMPL. The implementation required reverse engineering AMPL's syntax and developing a custom text-based parser. Such feats of software engineering are not needed to develop extensions to JuMP.

\subsection{Extension for robust optimization}
Robust optimization (RO) is a methodology for addressing uncertainty in optimization problems that has grown in popularity over the last decade (for a survey, see~\cite{ROSurvey}). The RO approach to uncertainty models the \textit{uncertain parameters} in a problem as belonging to an \textit{uncertainty set}, instead of modeling them as being drawn from probability distributions. We solve an RO problem with respect to the worst-case realization of those uncertain parameters over their uncertainty set, i.e.
\begin{alignat}{2}
\label{genericro}
\min_{x\in X} \quad & f \left( x \right) \\
\text{subject to} \quad &
g \left( x,\xi \right) \leq 0 \quad \forall \xi \in \Xi \nonumber
\end{alignat}
where $x$ are the decision variables, $\xi$ are the uncertain parameters drawn from the uncertainty set $\Xi$, $f : X \to \mathbb{R}$ is a function of $x$ and $g: X \times \Xi \to \mathbb{R}^k$ is a vector-valued function of both $x$ and $\xi$. Note that constraints which are not affected by uncertainty are captured by the set $X$. As the uncertainty set $\Xi$ is typically not a finite set of scenarios, RO problems have an infinite set of constraints. This is usually addressed by either reformulating the RO problem using duality to obtain a \textit{robust counterpart}, or by using a cutting-plane method that aims to add only the subset of constraints that are required at optimality to enforce feasibility~\cite{Cuts}.

JuMPeR~\cite{iainthesis} is an extension for JuMP that enables modeling RO problems directly by introducing the \texttt{Uncertain} modeling primitive for uncertain parameters. The syntax is essentially unchanged from JuMP, except that constraints containing only \texttt{Uncertain}s and constants are treated distinctly from other constraints as they are used to define the uncertainty set. JuMPeR is then able to solve the problem by either reformulation or the cutting-plane method, allowing the user to switch between the two at will. This is an improvement over both directly modeling the robust counterpart to the RO problem and implementing a cutting-plane method, as it allows users to experiment with different uncertainty sets and solution techniques with minimal changes to their code. 
Building JuMPeR on top of JuMP makes it more useful than a dedicated RO modeling tool like ROME~\cite{ROME} as users can smoothly transition from a deterministic model to an uncertain model and can take advantage of the infrastructure developed for JuMP to utilize a wide variety of solvers.
It also benefits from JuMP's efficient model construction, offering some performance advantages over YALMIP's robust modeling capabilities~\cite{lofberg2012automatic}.

\subsection{Extension for chance constraints}

For the final extension, consider chance constraints of the form
\begin{equation}\label{eq:chance}
\mathbb{P}(\xi^Tx \leq b) \geq 1-\epsilon,
\end{equation}
where $x$ is a decision variable and $\xi$ is a random variable. That is, $x$ is feasible if and only if the random variable $\xi^Tx$ is less than $b$ with high probability. Depending on the distribution of $\xi$, the constraint may be intractable and nonconvex; however, for the special case of $\xi$ jointly Gaussian with mean $\mu$ and covariance matrix $\Sigma$, it is convex and representable by conic-quadratic inequalities. Bienstock et al.~\cite{ccopf-sirev} observed that it can be advantageous to implement a custom cutting-plane algorithm similar to the case of robust optimization. The authors in~\cite{ccopf-sirev} also examined a more conservative \textit{distributionally robust} model where we enforce that~\eqref{eq:chance} holds for a family of Gaussian distributions where the parameters fall in some uncertainty set $\mu \in U_\mu, \Sigma \in U_\Sigma$.

JuMPChance is an extension for JuMP which provides a natural algebraic syntax to model such chance constraints, hiding the algorithmic details of the chance constraints from users who may be practitioners or experts in other domains. Users may declare Gaussian random variables and use them within constraints, providing $\epsilon$ though a special \texttt{with\_probability} parameter. JuMPChance was used to evaluate the distributionally robust model in the context of optimal power flow under uncertainty from wind generation, finding that the increased conservatism may actually result in realized cost savings given the inaccuracy of the assumption of Gaussianity~\cite{JuMPChanceCaseStudy}.

\section{Interactivity and visualization}\label{sec:ijulia}
Although we have focused thus far on efficiently and intuitively communicating optimization problems to a solver, equally as important
is a convenient way to interpret, understand, and communicate the solutions obtained. For many use cases, Microsoft Excel and similar spreadsheet systems provide a surprisingly versatile environment for optimization modeling~\cite{SolverStudio}; one reason for their continuing success is that it is trivial to interactively manipulate the input to a problem and visualize the results, completely within the spreadsheet. Standalone commercial modeling systems, while providing a much better environment for handling larger-scale inputs and models, have in our opinion never achieved such seamless interactivity. Notably, however, AIMMS~\cite{aimms}, a commercial AML, enables users to create interactive graphical user interfaces. We highlight, however, that both AIMMS and Excel-based solutions like SolverStudio~\cite{SolverStudio} require commercial software and are available only for the Windows operating system.

Many in the scientific community are beginning to embrace the ``notebook'' format for both research and teaching~\cite{IPythonNature}.  Notebooks allow users to mix code, rich text, \LaTeX\ equations, visualizations, and interactive widgets all in one shareable document, creating compelling narratives which do not require any low-level coding to develop. Jupyter~\cite{Jupyter}, in particular, contains the IJulia notebook environment for Julia and therefore JuMP as well. Taking advantage of the previously demonstrated speed of JuMP, one can easily create notebooks that embed large-scale optimization problems, which we will illustrate with two examples in this section. We believe that notebooks provide a satisfying solution in many contexts to the longstanding challenge of providing an interactive interface for optimization.

\subsection{Example: Portfolio Optimization}

One of the classic problems in financial optimization is the Markowitz portfolio optimization problem~\cite{Markowitz} where we seek to optimally allocate funds between $n$ assets. The problem considers the mean and variance of the return of the resulting portfolio, and seeks to find the portfolio that minimizes variance such that mean return is at least some minimal value. This is a quadratic optimization problem with linear constraints. It is natural that we would want to explore how the optimal portfolio's variance changes as we change the minimum return: the so-called \textit{efficient frontier}.

\begin{figure}
\centering
\includegraphics[width=\textwidth]{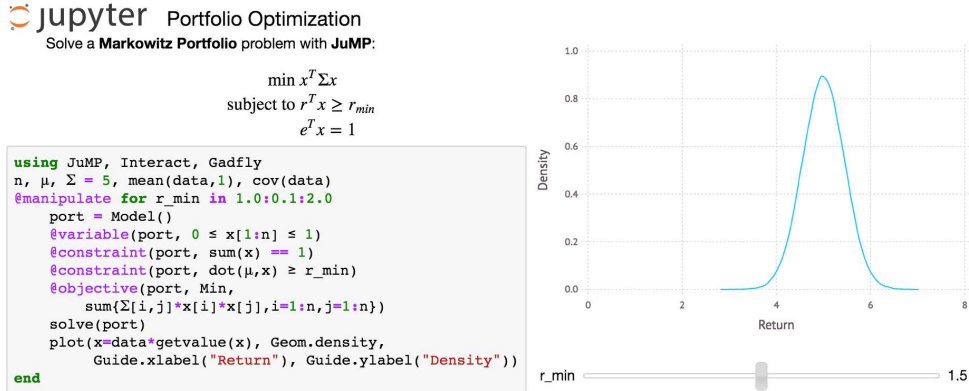}

\caption{A Jupyter (IJulia) Notebook for a Markowitz portfolio problem~\cite{Markowitz} that combines rich text with equations, Julia/JuMP code, an interactive widget, and a visualization. Moving the $r_{min}$ slider re-solves the optimization problem to find a new portfolio, and the plot is updated to show the historical distribution of returns that would have been obtained with the portfolio.}
\label{fig:ijuliaport}
\end{figure}

In Figure~\ref{fig:ijuliaport} we have displayed a small notebook that solves the Markowitz portfolio optimization problem. The notebook begins with rich text describing the formulation, after which we use JuMP to succinctly express the optimization problem. The data is generated synthetically, but could be acquired from a database, spreadsheets, or even directly from the Internet. The Julia package Interact.jl~\cite{InteractJL} provides the \texttt{@manipulate} syntax, which automatically generates the minimum return slider from the definition of the \texttt{for} loop. As the user drags the slider, the model is rebuilt with the new parameter and re-solved, enabling easy, interactive experimentation. The visualization (implemented with the Gadfly~\cite{GadflyJL} package) of the distribution of historical returns that would have been obtained with this optimal portfolio is also regenerated as the slider is dragged.

\subsection{Example: Rocket Control}

A natural goal in aerospace engineering is to maximize the altitude attained by a rocket in flight. This problem was possibly first stated by Goddard~\cite{goddardrocket}, and has since become a standard problem in control theory, e.g. \cite{brysondynopt}. The ``Goddard Rocket'' optimization problem, as expressed in~\cite{COPS3}, has three state variables (altitude, velocity, and remaining mass) and one control (thrust). The rocket is affected by aerodynamic drag and gravity, and the constraints of the problem implement the equations of motion (discretized by using the trapezoidal rule).

\begin{figure}
\centering
\includegraphics[width=0.8\textwidth]{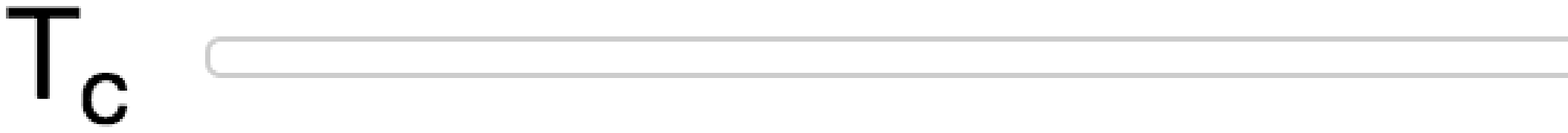}
\vspace{0.05cm}
\includegraphics[width=\textwidth]{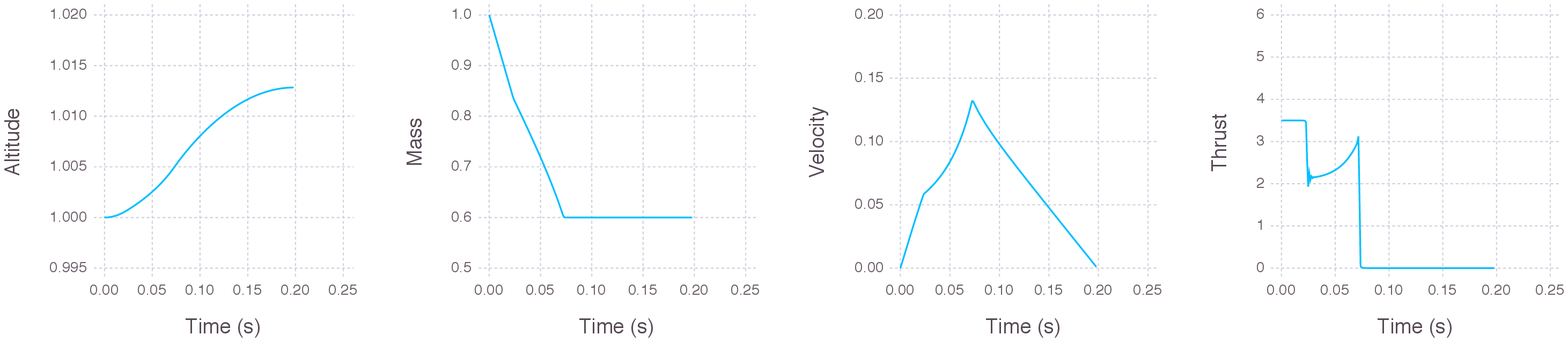}

\vspace{0.2cm}

\includegraphics[width=0.8\textwidth]{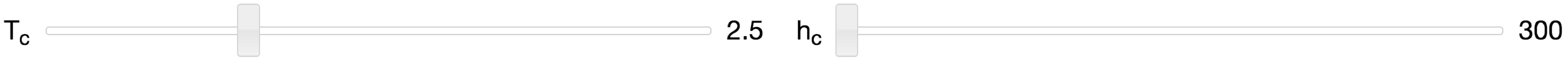}
\vspace{0.05cm}
\includegraphics[width=\textwidth]{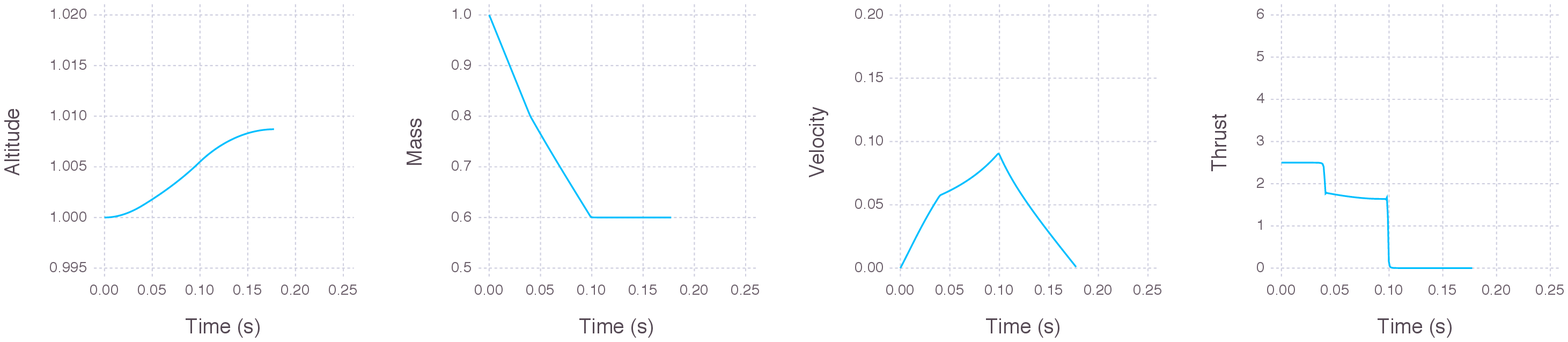}
\caption{Visualization of the states (altitude,mass,velocity) and the control (thrust) for a rocket optimal control problem. The top set of figures is obtained for the thrust and drag parameters, resp., ${T_c=3.5}$ and ${h_c=500}$, and the second are obtained for the parameters ${T_c=2.5, h_c=300}$, with all units normalized and dimensionless. We can see that the increased drag and reduced maximum thrust in the bottom set of figures has a substantial impact on maximum altitude and leads to a very different thrust profile.}
\label{fig:ijuliarocket}
\end{figure}

We have implemented the optimization problem with JuMP in an IJulia notebook. Moreover we have used Interact.jl to allow the user to explore the effects of varying the maximum thrust (via $T_c$) and the coefficient that controls the relationship between altitude and drag ($h_c$). The JuMP code is omitted for the sake of brevity, but the sliders and plots of the state and control over time are displayed in Figure~\ref{fig:ijuliarocket}. The model is re-solved with the new parameters every time the user moves the sliders; this takes about a twentieth of a second on a laptop computer, enabling real-time interactive exploration of this complex nonlinear optimization model.

\section*{Supplementary materials}
The benchmark instances used in Sections~\ref{sec:linear} and~\ref{section:nonlinear} and the notebooks presented in Section~\ref{sec:ijulia} are available as supplementary materials at \url{https://github.com/mlubin/JuMPSupplement}. The site \url{http://www.juliaopt.org/} is the homepage for JuMP and other optimization-related projects in Julia.

\section*{Acknowledgements}
We thank Paul Hovland and Jean Utke for helpful discussions on automatic differentiation. We thank Juan Pablo Vielma, Chris Coey, Chris Maes, Victor Zverovich, and the anonymous referees for their comments on this manuscript which improved its presentation. This work would not be possible without the supportive community of Julia developers and users who are too many to name. We thank Carlo Baldassi, Jack Dunn, Joaquim Dias Garcia, Jenny Hong, Steven G. Johnson, Tony Kelman, Dahua Lin, Karanveer Mohan, Yee Sian Ng, Elliot Saba, Jo\~{a}o Felipe Santos, Robert Schwarz, Felipe Serrano, Madeleine Udell, Ulf Wors\o{}e (of Mosek), and David Zeng for significant contributions to solver interfaces in Julia. We thank Jarrett Revels for his work on the ForwardDiff.jl package, and Steven Dirkse for his help with the GAMS version of the minimum cost flow example.
Finally, we thank the many students in the Operations Research Center at MIT who have been early adopters of JuMP.
This material is based upon work supported by the National Science Foundation Graduate Research Fellowship under Grant No. 1122374.
M. Lubin was supported by the DOE Computational Science Graduate Fellowship, which is provided under grant number DE-FG02-97ER25308.

\bibliographystyle{siam}

\section{Appendix: Benchmark models for Section~\ref{sec:conicbenchmark}}

\subsection{\texttt{lqcp}}

The linear-quadratic control problem is Equation~(5.2-I) from \cite{mittelmann2001sufficient}. This model has a quadratic objective and linear constraints, and can be scaled by increasing the discretization (parameters $m$ and $n$) of the two-dimensional problem domain. For the purposes of benchmarking we measured the model generation time across a range of sizes, fixing $m=n$ and varying $n\in\left\{500,1000,1500,2000\right\}$.
In the notation below, we define $I=\left\{0,\dots,m\right\}$ to be the index set along the first dimension and $J=\left\{0,\dots,n\right\}$ as the index set for the second. We additionally define $I^\prime \leftarrow I \setminus\left\{ m\right\}$ and $J^\prime \leftarrow J \setminus\left\{0,n\right\}$, with all other parameters defined as in the above reference.

\begin{alignat*}{2}
\min_{\mathbf{u},\mathbf{y}}\quad & \frac{1}{4}\Delta_{x}\left(\left(y_{m,0}-y_{0}^{t}\right)^{2}+2\sum_{j=1}^{n-1} \left(y_{m,j}-y_{j}^{t}\right)^{2}+\left(y_{m,n}-y_{n}^{t}\right)^{2}\right)+ & \quad\\
 & \qquad\frac{1}{4}a\Delta_{t}\left(2\sum_{i=1}^{m-1} u_{i}^{2}+u_{m}^{2}\right)\\
\text{s.t.}\quad & \nicefrac{1}{\Delta_{t}}\left(y_{i+1,j}-y_{i,j}\right)=\\
 & \qquad\frac{1}{2h_{2}}\left(y_{i,j-1}-2y_{i,j}+y_{i,j+1}+y_{i+1,j-1}-2y_{i+1,j}+y_{i+1,j+1}\right) &  & \forall i\in I^\prime,j\in J^\prime\\
 & y_{0,j}=0 &  & \forall j\in J\\
 & y_{i,2}-4y_{i,1}+3y_{i,0}=0 &  & \forall i\in I\\
 & \nicefrac{1}{2\Delta_{x}}\left(y_{i,n-2}-4y_{i,n-1}+3y_{i,n}\right)=u_{i}-y_{i,n} &  & \forall i\in I\\
 & -1\leq u_{i}\leq1 &  & \forall i\in I\\
 & 0\leq y_{i,j}\leq1 &  & \forall i\in I,j\in J
\end{alignat*}

\subsection{\texttt{fac}}\label{sec:facility}

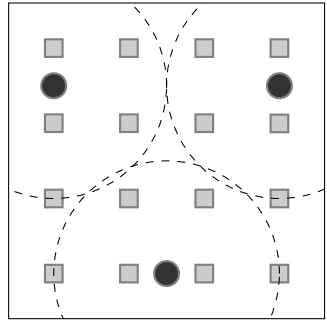
\begin{figure}
\centering

\begin{tikzpicture}
[scale=3,
 customer/.style={rectangle,draw=black!50,fill=black!20,thick},
 facility/.style={circle,draw=black!50,fill=black!80,thick}]
\node [customer] at (0.000,0.000) {};
\node [customer] at (0.333,0.000) {};
\node [customer] at (0.667,0.000) {};
\node [customer] at (1.000,0.000) {};
\node [customer] at (0.000,0.333) {};
\node [customer] at (0.333,0.333) {};
\node [customer] at (0.667,0.333) {};
\node [customer] at (1.000,0.333) {};
\node [customer] at (0.000,0.667) {};
\node [customer] at (0.333,0.667) {};
\node [customer] at (0.667,0.667) {};
\node [customer] at (1.000,0.666) {};
\node [customer] at (0.000,1.000) {};
\node [customer] at (0.333,1.000) {};
\node [customer] at (0.667,1.000) {};
\node [customer] at (1.000,1.000) {};

\node [facility] at (0.000,0.833) {};
\node [facility] at (1.000,0.833) {};
\node [facility] at (0.500,0.000) {};

\begin{scope}
\clip (-0.2,-0.2) rectangle ( 1.2, 1.2);
\draw [dashed] (0.000,0.833) circle (0.500);
\draw [dashed] (1.000,0.833) circle (0.500);
\draw [dashed] (0.500,0.000) circle (0.500);
\end{scope}

\draw (-0.2,-0.2) rectangle ( 1.2, 1.2);

\end{tikzpicture}

\caption{One possible optimal solution to the facility location problem with a four-by-four grid of customers (rectangles) and three facilities (circles). The dotted circles show the maximum distance between any customer and it's closest facility, which is the objective.}
\label{fig:facloc}
\end{figure}

The \texttt{fac} problem is a variant on the classic facility location problem~\cite{owen1998strategic}: given customers 
(indexed by $c \in \left\{ 1,\dots,C \right\}$)
located at the points
$x_{c} \in \mathbb{R}^K$, 
locate facilities
(indexed by $f \in \left\{ 1,\dots,F \right\}$)
at the points $y_{f} \in \mathbb{R}^K$ 
such that the maximum distance between a customer and its nearest facility is minimized. This problem can be expressed most naturally in the form of a mixed-integer second-order cone problem (MISOCP), and a solved example of this problem is presented in Figure~\ref{fig:facloc}. We generated the problem data deterministically to enable fair comparison across the different languages: the customers are placed on a two-dimensional grid ($K=2$)
$i \in \left\{ 0,\dots,G \right\}$ by 
$j \in \left\{ 0,\dots,G \right\}$, 
with the points $x_c$ spaced evenly over the unit square $\left[0,1\right]^2$.
The problem size is thus parametrized by the grid size $G$ and the number of facilities $F$, with the number of variables and constraints growing proportional to $F \cdot G^2$. For the purposes of benchmarking we measured the model generation with fixed $F=G$ and varying $F\in\left\{25,50,75,100\right\}$. 

\begin{alignat}{2}
\label{benchmodelsocp}
\min_{d,\mathbf{y},\mathbf{z}} \quad & d \\
\text{subject to} \quad & 
d \geq \left\Vert x_c - y_f \right\Vert_2 -
    M \left( 1 - z_{c,f} \right) & \quad &\forall c,f \nonumber \\
 & \sum_{f=1}^F z_{c,f} = 1 & \quad & \forall c \nonumber \\
 & z_{c,f} \in \left \{ 0,1 \right \} & \quad & \forall c,f, \nonumber
\end{alignat}
where 
\begin{equation*}
M = \max_{c,c^\prime} \left\Vert x_{c} - x_{c^\prime} \right\Vert_2
\end{equation*}
and $z_{c,f}$ is a binary indicator variable that is 1 if facility $f$ is closer to customer $c$ than any other facility, 0 otherwise. As needed, we translate the conic-quadratic constraint in~\eqref{benchmodelsocp} to an equivalent quadratic form depending on the modeling system.

\section{Benchmark models for Section~\ref{sec:nlbenchmark}}

\subsection{\texttt{clnlbeam}}

The first model, \texttt{clnlbeam}, is a nonlinear beam control problem obtained from Hans Mittelmann's AMPL-NLP benchmark set\footnote{\url{http://plato.asu.edu/ftp/ampl-nlp.html} accessed July 7, 2016.}; see also~\cite{clnlbeam}. It can be scaled by increasing the discretization of the one-dimensional domain through the parameter $n$. We test with $n \in \{5000,50000,500000\}$. The model has $3n$ variables, $2n$ constraints, and diagonal Hessians. The algebraic representation follows below.

\begin{align*}
	\min_{t,x,u \in \mathbb{R}^{n+1}} \quad & \sum_{i=1}^n \left[\frac{h}{2}(u_{i+1}^2 + u_i^2) + \frac{\alpha h}{2}(\cos(t_{i+1}) + \cos(t_i))\right]\\
	\text{subject to}\quad& x_{i+1} - x_{i} - \frac{1}{2n}(\sin(t_{i+1}) + \sin(t_i)) = 0\quad i = 1,\ldots,n\\
		& t_{i+1} - t_i - \frac{1}{2n}u_{i+1} - \frac{1}{2n}u_i = 0\quad i=1,\ldots,n\\
	 &-1 \leq t_i \leq 1,\quad -0.05 \leq x_i \leq 0.05\quad i = 1,\ldots,n+1\\
	 &x_1 = x_{n+1} = t_1 = t_{n+1} = 0.
\end{align*}

\subsection{\texttt{acpower}}

The second model is a nonlinear AC power flow model published in AMPL format by Artelys Knitro\footnote{\url{https://web.archive.org/web/20150105161742/http://www.ziena.com/elecpower.htm} accessed July 7, 2016.}. The objective is to minimize active power losses
\begin{equation}
\sum_k \left[g_k + \sum_{m} V_kV_m (G_{km}\cos(\theta_k-\theta_m) + B_{km}\sin(\theta_k-\theta_m))\right]^2
\end{equation}
subject to balancing both active and reactive power loads and demands at each node in the grid, where power flows are constrained by the highly nonlinear Kirchoff's laws. The parameter $g_k$ is the active power load (demand) at node $k$, $V_k$ is the voltage magnitude at node $k$, $\theta_k$ is the phase angle, and $Y_{km} = G_{km} + iB_{km}$ is the complex-valued admittance between nodes $k$ and $m$, which itself is a complicated nonlinear function of the decision variables.  Depending on the physical characteristics of the grid, some values (e.g., $V_k$) may be decision variables at some nodes and fixed at others. This model is quite challenging because of the combination of nonlinearity and network structure, which yields a highly structured Hessian.

We translated the AMPL model provided by Artelys Knitro to JuMP, GAMS, and Pyomo. The base instance has a network with 662 nodes and 1017 edges; there are 1489 decision variables, 1324 constraints, and the Hessian (of the Lagrangian) has 8121 nonzero elements.  We artificially enlarge the instances by duplicating the network 10-fold and 100-fold, which results in proportional increases in the problem dimensions.

\end{document}